\documentclass[twoside, 12pt]{article}
\usepackage[noadjust,sort]{cite}
\usepackage{amsfonts}
\usepackage{amsmath}
\usepackage{amssymb}
\usepackage{blkarray}
\usepackage{multirow}
\usepackage{graphics}

\usepackage{amsthm}
\usepackage{varioref}
\usepackage[colorlinks=true,allcolors=blue]{hyperref}

\usepackage[nameinlink,capitalize,noabbrev]{cleveref}
\crefname{equation}{}{}

\usepackage{pgf,tikz}
\usepackage{mathrsfs}
\usetikzlibrary{arrows}
\definecolor{wrwrwr}{rgb}{0.3803921568627451,0.3803921568627451,0.3803921568627451}
\definecolor{rvwvcq}{rgb}{0.08235294117647059,0.396078431372549,0.7529411764705882}
\newcommand{\R}{\mathbb{R}}
\newtheorem{theorem}{Theorem}[section]
\newtheorem{remark}[theorem]{Remark}
\newtheorem{example}[theorem]{Example}
\newtheorem{lemma}[theorem]{Lemma}
\newtheorem{corollary}[theorem]{Corollary}
\newtheorem{definition}[theorem]{Definition}
\newtheorem{proposition}[theorem]{Proposition}

\newtheorem{conjecture}{Conjecture}

\newcommand{\bea}{\begin{eqnarray}}
\newcommand{\eea}{\end{eqnarray}}

\newcommand{\comment}[1]{}
\def \bpm{\begin{pmatrix}}
\def \epm{\end{pmatrix}}
\def \bd{\begin{definition}}
\def \ed{\end{definition}}
\def \bcc{\begin{conjecture}}
\def \ecc{\end{conjecture}}
\def \bt{\begin{theorem}}
\def \et{\end{theorem}}
\def \bl{\begin{lemma}}
\def \el{\end{lemma}}
\def \bc{\begin{corollary}}
\def \ec{\end{corollary}}
\def\be#1\ee{\begin{align}#1\end{align}}
\def\beq #1\eeq {\begin{align*}#1\end{align*}}

\def \ben{\begin{enumerate}}
\def \een{\end{enumerate}}
\def \ba{\begin{array}}
\def \ea{\end{array}}
\def \bp{\begin{proposition}}
\def \ep{\end{proposition}}
\def \bx{\begin{example}}
\def \ex{\end{example}}
\def \br{\begin{remark}}
\def \er{\end{remark}}
\def \bdsc{\begin{description}}
\def \edsc{\end{description}}
\def\pf{{\it \bf Proof. }}
\def \qed {\hfill \vrule height6pt width6pt depth0pt}
\def\hs{\hspace{.3cm}}
\def\vs{\vskip .3cm}
\def\ds{\displaystyle}
\def\1{1\!\!1}
\def\J{\mathbb{J}}
\def\x{\mathbf{x}}
\def\r{\rho}
\def\D{\mathcal{D}}
\title{On the largest two and smallest six distance Pareto eigenvalues of a graph.}
\author{Deepak Sarma \\Department of Mathematical Sciences, \\ Tezpur University, Tezpur-784028, India.\\
Email address: \url{deepaks@tezu.ernet.in}}

\date{}

\begin{document}
\maketitle
\begin{abstract}
In this article, we establish some bounds involving the largest two distance Pareto eigenvalues of a connected graph. Also we characterize all possible values for smallest six distance Pareto eigenvalues of a connected graph. \\


% American Mathematical Society Subject Classification
\end{abstract}

\maketitle
\bigskip

\section{Introduction}

\hs All our graphs are finite, undirected, connected and simple.
Let $G$ be a  graph on vertices $\{1,2,\ldots,n\}$. At times, we use $V(G)$ and
$E(G)$ to denote the set of vertices and the set of edges of $G$, respectively. For $i, j\in V(G)$,
the {\em distance} between $i$ and $j$, denoted by $d_G(i,j)$ or simply $d_{ij}$, is the length of
a shortest path from $i$ to $j$ in $G$. The {\em distance matrix} of $G$, denoted by
$\mathcal{D}(G)$ is the $n\times n$ matrix with $(i,j)$-th entry $d_{ij}$.

For a column vector $x=(x_1,\ldots,x_n)^T\in\mathbb{R}^n$ we have
\be\label{quadratic-form}
x^T\D(G)x=\sum_{1\le i<j\le n}d_{ij}x_ix_j. \ee

If vertices $i$ and $j$ are adjacent, we write $i\sim j.$ Degree of a vertex $v$ in a graph $G$ will be denoted by $d_G(v).$ By pendent vertex of a graph we mean a vertex of degree 1. The transmission, denoted by $Tr(v)$ of a vertex $v$ is the sum of the distances from $v$ to all other vertices in $G$. The diameter of a connected graph $G$ denoted by $diam(G)$ is the maximum distance between any two vertices in $G,$ i.e. $diam(G)$ is the largest entry of $D(G).$ A \textit{clique} of a graph is a maximal complete subgraph and \textit{clique number} of a graph is the order of a maximal clique. We denote \textit{clique number} of a graph $G$ by $\omega(G).$
By $K_n,P_n,C_n,S_n,~ \hbox{and } W_n$ we respectively mean the Complete graph, Path graph, Cycle graph, Star graph and Wheel graph with $n$ vertices. The complete bipartite graph with bipartition size $m$ and $n$ is represented by $K_{m,n}.$ By $S_n^+$ we represent the graph obtained by adding an edge between any two independent vertices in $S_n$. If $V_1\subseteq V(G)$ and $E_1\subseteq E(G),$ then by $G-V_1$ and $G-E_1$ we mean the graphs obtained from $G$ by deleting the vertices in $V_1$ and the edges $E_1$ respectively. In particular case when $V_1=\{u\}$ or $E_1=\{e\},$ we simply write $G-V_1$ by $G-u$ and $G-E_1$ by $G-e$ respectively. By $K_n-e$ is the graph obtained from $K_n$ by removing any one edge of it. The graph obtained from $G$ and $H$ by identifying $u\in G$ and $v\in H$ is denoted by $G_u*H_v$. When there is no confusion of vertices we write $G*H$ for the coalescence of the graphs $G$ and $H.$
%If $v$ is a vertex of a tree $T,$ then the components of $T-v$ are called the \textit{branches} of $T$ at $v.$

By spectral radius of a symmetric matrix $M$, we mean its largest eigenvalue and denote it by $\rho(M).$
Note that for a connected graph $G$, $\D(G)$ is irreducible nonnegative
matrix. Thus by the Perron-Frobenius Theorem, $\rho(\D)$ is simple, and there is a
positive eigenvector of $\D(G)$ corresponding to $\rho(\D)$.
Such eigenvectors corresponding to $\rho(\D)$ is called {\em Perron vector}
of $\D(G)$. By an eigenvector we mean a unit eigenvector and by $\mathbb{M}_n,$ we denote the class of all real matrices of order n. We use the notation $A\geq 0$ to indicate that each component of the matrix $A$ is nonnegative. Furthermore in places we write $A\geq B$ to mean $A-B\geq 0.$

\bd
A real number $\lambda$ is said to be a Pareto eigenvalue of $A\in \mathbb{M}_n$ if there exists a nonzero vector $\x(\geq 0)\in \R^n$ such that
\beq A\x \geq\lambda \x \quad and \quad \lambda=\frac{\x^TA\x}{\x^T\x}, \eeq
also we call $\x$ to be a Pareto eigenvector of $A$ associated with Pareto eigenvalue $\lambda$.
\ed
A Pareto eigenvalue of $\D(G)$ of a graph $G$ will be called as distance Pareto eigenvalue of $G.$ Fernandes at.el. in \cite{fjt17} and Seeger at \cite{see18} studied the Pareto eigenvalues of adjacency matrix of a graph. Pareto eigenvalue of the distance matrix of a connected graph was first studied in \cite{ns18}. In this article we study something more about distance Pareto eigenvalues.
\par This article is organized as follows. Some basic results of distance Pareto eigenvalues of a graph are discussed in \cref{prelim}. We establish some bounds of the difference and ratio of the largest two distance Pareto eigenvalues of a graph in \cref{diff and ratio}.  We characterize all possible smallest five distance Pareto eigenvalues of a connected graph in \cref{upto 5th}. Finally in \cref{6th}, we find the possible values of sixth smallest distance Pareto eigenvalue of a connected graph with at least $5$ vertices.

\section{Preliminaries and basic results}\label{prelim}

For a square matrix $A$, we use the symbol $A(i)$ for the principal submatrix of $A$ obtained by deleting $i-th$ row and column  of $A.$
In particular if $\D(G)$ be the distance matrix of a graph $G$ then by $\D(i)$ we will denote the principal submatrix of $\D$ obtained by deleting row and column corresponding to vertex $i$ of $G.$ By $\1$ we denote the column vector of all ones and by $\J$ the matrix of all ones of appropriate size.\\

\par For a matrix $A$ we use $\r_k(A)$ and $\mu_k(A)$ to denote the k-th largest and k-th smallest Pareto eigenvalue of $A.$
 For a connected graph $G$ we simply write $\r_k(G)$ and $\mu_k(G)$ to mean $\r_k(\D(G))$ and $\mu_k(\D(G))$ respectively.

 \bl{ [Weyl's Inequalities]}\cite{hj05}\label{w} Let $\lambda_i(M)$ denote the $i$-th largest eigenvalue of a real symmetric matrix $M.$ If $A$ and $B$ are two real symmetric matrices of order $n,$
then $$\lambda_1(A)+\lambda_i(B)\ge \lambda_i(A+B)\ge \lambda_n(A)+\lambda_i(B) \mbox{ for } i=1,2,\ldots, n.$$\el

\bl \cite{min88}\label{mi} If $A$ is an irreducible matrix and $A\geq B\geq 0, A\neq B,$ then $\r(A)>\r(B).$
\el

\bl\cite{hj05}\label{r} If $A$ is a symmetric $n\times n$ matrix with $\lambda_1$ as the largest eigenvalue then for any normalized vector
$\x\in \R^n (\x\ne 0)$,
\beq \x^t A \x\le \lambda_1.\eeq The equality holds if and only if $\x$ is an eigenvector of $A$ corresponding to the eigenvalue $\lambda_1.$\el

Putting $\x=\frac{1}{\sqrt{n}}\1$ in \cref{r} we get the following result as a Corollary.

\bc \label{avg} If $A$ is a symmetric $n\times n$ matrix, then
\be \label{av} \r(A)\ge \bar{R}, \ee \\
 where $\bar{R}$ is the average row sum of the matrix $A.$ The equality in \cref{av} holds if and only if all the row sums of $A$ are equal.
\ec

\bt \cite{sv11}\label{imp}
The scalar $\lambda \in \R$ is a Pareto eigenvalue of $A\in \mathbb{M}_n$ if and only if there exist a nonempty set $J\subset \{ 1, 2, \ldots, n\}$ and a vector $\xi\in \R^{|J|}$ such that

%\begin{align}
\beq
A^J\xi &=\lambda\xi, \\
\xi_j  &>0 \quad \forall j \in J.  \nonumber \\
 \sum_{j\in J}a_{i,j}\xi_j & \geq 0 \quad \forall i\notin J. \nonumber
\eeq
%\end{align}
 Furthermore, a Pareto eigenvector $\x$ associated to $\lambda$ is constructed by setting
%\begin{align*}
\be
\x_j =
  \begin{cases}
   \xi_j & \text{if } j\in J \nonumber \\
   0       & otherwise.
  \end{cases}
%\end{align*}
\ee
\et

\par From \cref{imp}, we get the following result similar to that of \cite[Theorem 1]{see18}.
\bt\cite{ns18} \label{main} The distance Pareto eigenvalues of a connected graph $G$ are given by

$\Pi(G)=\{ \r(A): A\in M\},$
where $M$ is the class of all principal sub-matrices of $\D(G).$
\et

\bl \label{k} \cite{ns18} For any positive integer $n,$ $\Pi(K_n)=\{0, 1, \ldots, n-1\}.$
\el

\bl \cite{ns18} \label{r2def} If $G$ is a connected graph with at least two vertices, then
\beq \r_2(G)=\max \{ \r(A): A\in P\},
\eeq

where $P=\{ (\D(G))(v): v\in V(G), d_v>1 \}$
\el

\bd
If $A$ and $B$ are two nonnegative matrices then we say that $A$ dominates $B$ if either of the following two cases hold
\ben
\item $A$ and $B$ are of same size and upto permutation similarity $A\geq B,$ $A\neq B.$
\item $A$ is permutation similar to $\left(\begin{array}{rr}
B & C \\
D & E
\end{array}\right)$ and at least one of C,D and E is a nonzero matrix.
\een
\ed

\bl \cite{ns18} \label{dom}
If $A$ and $B$ are two symmetric nonnegative irreducible matrices, then \\
$A$ dominates $B$ implies $\rho(A)>\rho(B).$
\el

\bl \cite{ns18} \label{Sn} There are $2(n-1)$ distance Pareto eigenvalues of $S_n$ and they are
\beq \mu_{2k}=2(k-1), \, \mu_{2k-1}=k-1+\sqrt{k^2-3k+3}
\text{  where }  k=1,\ldots, n-1.
\eeq
\el

\bt $\r_1(S_n)-\r_2(S_n) $ is a decreasing function in $n.$ \et
\pf Let \beq f(n)&=\r_1(S_n)-\r_2(S_n) \\
&= \sqrt{n^2-3n+3}-n+2. \qquad \Big[ \hbox{Using } \cref{Sn}   \Big]     \eeq
Then \beq f'(n)&=\frac{n-\frac{3}{2}-\sqrt{(n-\frac{3}{2})^2+\frac{3}{4}}   }{\sqrt{n^2-3n+3}}\\
& < 0 \quad \forall ~ n \in \mathbb{N} \eeq
This completes the proof. \qed

\bt $\r_2(S_n^+) =\frac{2n-7+\sqrt{(2n-1)^2-16}}{2} $\et
\pf Let $A=\D(S_n^+)(v),$ where $v$ is the vertex of $(S_n^+)$ of degree $n-1.$ Then upto permutation similarity, we get
\beq
A=\bpm  2(\J-I)_{n-3} & 2\J_{n-3,\, 2}\\
2\J_{2,\, n-3} & I_2 \epm
\eeq
Let $\x$ be the perron vector of $A.$ Then due to symmetry we can assume that all the components of $\x$ corresponding to pendent vertices equals $a$ and that of  the vertex of degree $n-1$ is $b.$ If $\r$ is the largest eigenvalue of $A,$ then from eigen equations we have

\be \label{Sn+1} 2(n-4)a+2\times2b &=\r a \\
\label{Sn+2} 2(n-3)a+b &=\r b
 \ee
Solving \cref{Sn+1,Sn+2} we get
\beq \r^2-(2n-7)\r +16-6n=0 \\
\therefore \quad \r =\frac{2n-7+\sqrt{(2n-5)(2n+3)}}{2} \eeq

Upto permutation similarity there are exactly 3 principal submatrices of $\D(S_n^+)$ of order $n-1.$ But $A$ dominates the other two, hence the result follows from \cref{r2def}. \qed

By direct varification we get the following result regarding the second largest distance pareto eigenvalue among all unicyclic graphs of order at most six.

\bt Among all unicyclic graph of order $n\le 6,$ $\r_2$ is minimum for $C_n$ and second minimum for $S_n^+.$ \et

\bcc Among all unicyclic graph $G$ of order $n\ge 7,$ $\r_2 \ge \frac{2n-7+\sqrt{(2n-1)^2-16}}{2} $ with equality if and only if $G \cong   S_n^+.$ \ecc

\bt $\r_2 (W_n)=2(n-3).$  \et
\pf Observe that upto permutation similarity there are only two distinct submatrix of $\D(W_n)$ of order $n-1.$ Let $A$ be the principal submatrix of $\D(W_n)$ of order $n-1$ obtained by removing row and column corresponding to vertex of degree $n-1$ and $B$ be any other principal submatrix of $\D(W_n)$ of order $n-1.$ Then $B$ has constant row sum equal to $2+2(n-4)=2(n-3).$
Therefore $\r(B)=2(n-3).$ Since $B$ dominates $A,$ hence the result follows. \qed

\bt \label{r1r2} If $k$ be the minimum positive component of the distance Pareto eigenvector of a graph $G$ corresponding to $\r_2,$ then for any real $t>0$
\beq  \r_1 \ge \frac{\r_2+2tk(n-1)}{1+t^2}  \eeq
with equality if and only if   $t=\frac{\sqrt{n-1}}{\r_1}=\frac{\r_1-\r_2}{\sqrt{n-1}}.$ \et

\pf Let $\D(G)= \bpm A & y \\ y^T & 0\epm $ so that $\r_2(G)=\r(A)$ and $z$ be the normalized vector with $z^T A z=\r_2.$ For $t>0,$ we set $x=\bpm z \\ t\epm .$ Then we have

\be \label{r1r21} \r_1 &\ge \frac{x^T \D x}{x^Tx} \\
& = \frac{z^TAz +2tz^Ty}{1+t^2} \nonumber \\
&= \frac{\r_2+2tz^Ty}{1+t^2} \nonumber \\
\label{r1r22}  &\ge \frac{\r_2+2\,t\, k\, trace(v)}{1+t^2} \\
\label{r1r23}  &\ge  \frac{\r_2+2tk (n-1)}{1+t^2}
\ee
Thus the first part is done.
 Now if the equality holds then equality \cref{r1r23} gives
 \be \label{r1r24}  trace(v)=n-1 \Rightarrow y=\1 . \ee
 Equality \cref{r1r22} gives
 \be \label{r1r25} z=\frac{\1}{\sqrt{n-1}}  \ee
 Again equality in \cref{r1r21} gives
 \be \label{r1r26}  \D x=\r_1 x \ee
 Using \crefrange{r1r24}{r1r26} we get
 \beq (\r_1-\r_2)\frac{\1}{\sqrt{n-1}}=t \, \1  \text{ and } \frac {n-1}{\sqrt{n-1}}=\r_1\, t. \\
 \text{ Therefore} \quad t=\frac{\sqrt{n-1}}{\r_1}=\frac{\r_1-\r_2}{\sqrt{n-1}}.
 \eeq
\qed

\bd We define by $\mathcal{G}_n$ the class of all connected graphs $g$ of order $n$ so that if $\r_2(g)=\r(A)$ where  $A\in \mathbb{M}_{n-1} $ is a principal submatrix of $\D(g)$ then $A$ has all row (column) sums equal.  \ed
{\bf Note:} A graph $g\in \mathcal{G}_n$ if and only if all non zero components of distance Pareto eigenvector of $g$ corresponding to $\r_2$ are $\frac{1}{\sqrt{n-1}}.$

\bt If $k$ be the minimum positive component of the distance Pareto eigenvector of a graph $G$ of order $n$ corresponding to $\r_2(G),$ then
\beq \r_1 \ge \frac{\r_2+\sqrt{\r_2^2+4(n-1)(2k\sqrt{n-1}-1)}  }{2}  \eeq
equality holds if and only if $G\in \mathcal{G}_n.$ \et
\pf Taking $t=\frac{\sqrt{n-1}}{\r_1}$ in \cref{r1r2} we get
\beq \r_1 & \ge \frac{ \r_2 +2\frac{\sqrt{n-1}}{\r_1}k(n-1) }{1+\frac{n-1}{\r_1^2}}\\
\hbox{i.e. } \quad \r_1 &\ge \frac{\r_2+\sqrt{\r_2^2+4(n-1)(2k\sqrt{n-1}-1)}  }{2} \eeq

Now from \cref{r1r2} equality holds in the above expression if and only if $k=\frac{1}{\sqrt{n-1}}$ i.e. if and only if $G\in \mathcal{G}_n.$ \qed

\bt If $k$ is the minimum positive component of the distance Pareto eigenvector of a graph $G$ corresponding to $\r_2(G),$ then
\beq 2\r_1-\r_2 \ge 2k(n-1), \eeq
equality holds if and only if $G=K_2.$
\et
\pf Taking $t=1$ in \cref{r1r2} we get
\beq \r_1 \ge \frac{\r_2+2k(n-1)}{2}\\
\hbox{i.e. } 2\r_1-\r_2 \ge 2k(n-1).\eeq
Now by \cref{r1r2} equality holds if and only if
\beq 1=\frac{\sqrt{n-1}}{\r_1}=\frac{\r_1-\r_2}{\sqrt{n-1}}\\
\hbox{i.e. } \r_1 =\sqrt{n-1} \text{ and } \r_2=0 
\eeq
  which is possible only if $G=K_2.$
\qed

\bd \label{def} A vertex $v$ of a connected graph $G$ of order $n$ is called pyramidal if $d_v =n-1$ and $G-v$ is connected and regular. Besides we call a connected graph $G$ to be
pyramidal if there exist at least one pyramidal vertex in it.
\ed

\bt If $G$ be a connected graph of order $n$ with $v \in V(G)$ such that $Tr(v)$ is minimum, then
\beq \r_1 \ge \frac{Tr(v) -1+\sqrt{(Tr(v)-1)^2+4(n-1) } }{2}  \eeq
with equality if and only if $G$ is pyramidal. \et

\pf Let $\x=(x_1,x_2, \ldots , x_n )$ be the Perron-vector with $x_i=\min x_k, x_j=\min \{x_k >x_i\}.$
From eigenequations, we have

\be \label{r1b1} \r_1x_i \ge T_i x_j \text{ and } \r_1x_j \ge x_i +(T_j-1)x_j\ee

Now \cref{r1b1} gives
\beq \r_1(\r_1-T_j+1)&\ge T_i \\
\Rightarrow \r_1^2-(T_j-1)\r_1-T_i &\ge 0 \\
\Rightarrow \quad \r_1 & \ge \frac{T_j-1+\sqrt{(T_j-1)^2+4T_i } }{2} \\
&\ge \frac{T_j-1+\sqrt{(T_j-1 )^2+4(n-1)}}{2}
\eeq
Thus the first part is done.\\
Now if the equalities hold, then considering all the above equalities we get
\beq v_i \sim v_j, \,  x_k=x_j \, \forall \, k\ne i \text{ and } d_i=n-1. \eeq
\beq \hbox { Therefore } \quad \r_1 x_j &=x_i+(T_j-1)x_j\\
\Rightarrow \quad T_k&=T_j\, \,\, \forall k\ne i \\
\Rightarrow T_k &= d_k +2(n-d_k-1) \, \forall k \ne i \\
& = 2(n-1)-d_k \\
\Rightarrow d_k&=d_j\,\, \forall i\ne k \text{ and } T_i=n-1.
\eeq
Hence the result follows.
\qed

\comment{
\bt If $G$ is a graph of order $n$ with $diam(G)=d>2,$ then
\beq |\Pi(G)|\geq n+2(d-1).\eeq
\et
\pf If we consider $A_i$ as in Theorem \ref{n} for $i=3,\ldots , n,$ then we find $n+d-1$ distance Pareto eigenvalues of $G.$ Now  let
$B_k=\bpm
0 & 1 & k \\
1 & 0 & k-1 \\
k & k-1 & 0
\epm.$
 Then for each $k=2,3,\ldots,d,$ $B_k$ is a principal sub-matrix of $\D(P_d).$ But minimum and maximum row sum of $B_k$ are $k$ and $k+1$ respectively.

 Therefore $\r(B_k)\in (k,k+1).$
 Besides minimum row sum of $A_3$ is at least $d+1$ and therefore,
 $d<\r(B_d)<\r(A_3).$ Thus we get $n+d-1+(d-1)=n+2(d-1)$ distance Pareto eigenvalues of $G$ as follows
 \beq 0<1<2<\r(B_2)<3<\ldots <\r(B_{d-1})<d<\r(B_d)<\r(A_3)\ldots <\r(A_n).\eeq
Hence the result.\qed }

\section{Difference and ratio of largest two distance Pareto eigenvalues}\label{diff and ratio}

From \cite{ns18} we see that the second largest distance Pareto eigenvalues
is always greater than the second largest distance eigenvalue, we now study the difference(ratio) of the largest two distance Pareto eigenvalues of a graph.

\bt For any positive integer $n,$ $\r_1(C_n)-\r_2(C_n)<\frac{1}{n-1}\lfloor \frac{n^2}{4} \rfloor .$
\et
\pf Since all row(column) sum of $\D(C_n)$ are equal to $\lfloor \frac{n^2}{4} \rfloor ,$ therefore
\be \label{411} \hs \r_1(C_n)=\bigg\lfloor \frac{n^2}{4} \bigg\rfloor . \ee

Now upto permutation similarity all the principal sub-matrix of $\D(C_n)$ of order $n-1$ are equal and average row sum of any such matrix is \beq \frac{n\lfloor \frac{n^2}{4} \rfloor-2\lfloor \frac{n^2}{4} \rfloor}{n-1} &= \frac{n-2}{n-1}\bigg\lfloor \frac{n^2}{4} \bigg\rfloor \\
&=\bigg\lfloor \frac{n^2}{4}\bigg \rfloor-\frac{1}{n-1}\bigg\lfloor \frac{n^2}{4}\bigg \rfloor .\eeq

Besides all the row(column) sum of principal sub-matrix of $\D(C_n)$ of order $n-1$ are not equal. Therefore using \cref{avg} we get
\be \label{412}
 \r_2(C_n)>\bigg\lfloor \frac{n^2}{4}\bigg \rfloor-\frac{1}{n-1}\bigg\lfloor \frac{n^2}{4}\bigg \rfloor .
\ee
 From \cref{411,412} we have

\beq
\r_1(C_n)-\r_2(C_n)&< \bigg\lfloor \frac{n^2}{4}\bigg\rfloor-\bigg\lfloor \frac{n^2}{4}\bigg \rfloor+\frac{1}{n-1}\bigg\lfloor \frac{n^2}{4}\bigg \rfloor \\
&=\frac{1}{n-1}\bigg\lfloor \frac{n^2}{4} \bigg\rfloor .\eeq
\qed

\bcc If $G$ is a connected graph of order $n,$ then $ \r_1(G)-\r_2(G)\leq \r_1(C_n)-\r_2(C_n),$ equality holds if and only if $G=C_n.$
\ecc

\bt For any positive integer $n,$ $\r_1(S_n)-\r_2(S_n)=\sqrt{n^2-3n+3}-n+2  .$
\et
\pf From \cref{Sn} we have \[ \r_1(S_n)=n-2+\sqrt{n^2-3n+3} \hs \text{and} \hs \r_2(S_n)=2(n-2).\] Hence the result follows easily. \qed

\bcc \label{con3} If $G$ is a connected graph of order $n,$ then \[ \r_1(G)-\r_2(G)\geq \sqrt{n^2-3n+3}-n+2 ,\]
\par equality holds if and only if $G=S_n.$
\ecc

\bcc Among all connected graphs of order $n$ the sum of $k$ largest distance Pareto eigenvalue is minimum for $S_n$ and maximum for $P_n.$
\ecc 

\bcc{} Among all bipartite graph the sum of $k$ largest distance Pareto eigenvalue is minimum for $K_{\lfloor \frac{n}{2}\rfloor \lceil \frac{n}{2}\rceil}$
\ecc

\bt For any connected graph $G,$ $\r_1(G)-\r_2(G)< {\ds \min_{u\in V(G)} \sqrt{\sum_{v\in V(G)}d_{uv}^2} } .$
\et
\pf Let \[ \D(G)= \bpm
0 & \x^T \\
\x & E
\epm , \hs
 M=\bpm
0 & 0 \\
0 & E
\epm  \quad \text{and }
N= \bpm
0 & \x^T \\
\x & 0
\epm .
\]
\par Then clearly $\D=M+N,$ where $\D=\D(G).$
Now by Weyl's inequality (\ref{w}) we have
\be \label{430} \r(\D) \leq \r(M) + \r(N),\ee equality holds if and only if there is a vector which is at a time eigenvector corresponding to $\r(\D)$ of $\D$, $\r(M)$ of $M$ and $\r(N)$ of $N.$
\par Now we have \be
\label{431}\r(\D)=\r_1(G),\\
 \r(M)\leq \r_2(G)\\
 \text{and }\quad  \r(N)=\sqrt{\x^T\x}\ee

\par Also there is a particular $u\in V(G)$ such that \be \label{432} \x_i=d_{ui}\ee

On using \cref{431,432}, inequality \cref{430} reduces to the form
 \be \label{433}
\r_1(G)\leq \r_2(G)+\sqrt{\sum_{v\in V(G)}d_{uv}^2}.
 \ee
Taking minimum of $\sum_{v\in V(G)}d_{uv}^2$ over all vertices of $G$ in \cref{433}, we get
\be \label{434}
\r_1(G)-\r_2(G)\leq {\ds \min_{u\in V(G)} \sqrt{\sum_{v\in V(G)}d_{uv}^2} }. \ee

If possible suppose the equality holds in \cref{434}. Then $\exists y\in \mathbb{C}^n$ such that

\be
\label{435} \D y=\r(\D)y \\
My=\r(M)y \\
\label{436} Ny=\r(N)y
\ee

From \cref{435,436} we have
\be
\label{437} \sum_{v\in V(G)}d_{uv}y_v=\r(\D)y_u \\
\label{438} \sum_{v\in V(G)}d_{uv}y_v=\r(N)y_u
\ee

Equations \cref{437,438} suggests that $\r(\D) =\r(N),$ which is a contradiction to the fact that $\D$ dominates $N.$ Therefore the equality in \cref{434} can never hold. Hence the result follows.
\qed

\bl \label{ratio}
If $G$ be a connected graph and $\x$ be the normalized perron vector of $\D(G)$ then for any $v \in V(G)$,
\beq
\frac{\r_1}{\r_2}\le \frac{1-x_v^2}{1-2x_v^2}
\eeq
equality holds if and only if $x_u=\frac{d_{uv}}{\sqrt{(\r_1-\r_2)(2\r_1-\r_2)}} \quad \forall \, u\neq v$
\el
\pf
Upto permutation similarity we can take
$\D=\bpm
A & y \\
y^T & 0
\epm$ and $\x = \bpm  z \\ x_v \epm ,$
where $\D=\D(G)$ and $y_u=d_{uv}.$

Now $\D\x=\r_1\x$ gives
\be
\label{r2} Az+x_vy &=\r_1z \\
\label{r3} y^Tz &=\r_1x_v
\ee
Also as $\x$ is normalized vector so we have
\be
\label{r4} z^Tz=1-x_v^2
\ee

From \cref{r2,r3,r4}, we get
\be
\label{r5} z^TAz=\r_1(1-2x_v^2)
\ee

But from \cref{r}, we have
\be
\label{r6} \frac{z^TAz}{z^Tz}\le \r_2
\ee
Equations \cref{r5} and \cref{r6} together gives
\be
\label{r7} \frac{\r_1}{\r_2}\le \frac{1-x_v^2}{1-2x_v^2}
\ee
Thus the first part is done.\\
\par Now suppose the equality holds in \cref{r7}, then equality must hold in \cref{r6} as well. Therefore we have
\beq  \quad Az=\r_2z. \eeq

Using equation \cref{r2} we get
\be
\label{r8} z=\frac{x_v}{\r_1-\r_2}y
\ee
But equality in \cref{r7} gives
\beq
x_v=\sqrt{\frac{\r_1-\r_2}{2\r_1-\r_2}}
\eeq
Therefore equation \cref{r8} reduces to
\beq
z=\frac{y}{\sqrt{(\r_1-\r_2)(2\r_1-\r_2)}}
\eeq
which is again equivalent to the form

\beq
x_u=\frac{d_{uv}}{\sqrt{(\r_1-\r_2)(2\r_1-\r_2)}} \quad \forall \, u\neq v
\eeq
 \par Conversely if we assume $x_u=\frac{d_{uv}}{\sqrt{(\r_1-\r_2)(2\r_1-\r_2)}} \quad \forall \, u\neq v,$
then from equations \cref{r3,r4}, we get
\be
\label{r9}  \frac{\sum_{u\in V(G)d_{uv}^2}}{\sqrt{(\r_1-\r_2)(2\r_1-\r_2)}}&=\r_1x_v \\
\text{ and } \quad \label{r10} \frac{\sum_{u\in V(G)d_{uv}^2}}{(\r_1-\r_2)(2\r_1-\r_2)}&=1-x_v^2
\ee
Simplifying equations \cref{r9,r10}, we get
$\frac{\r_1}{\r_2}= \frac{1-x_v^2}{1-2x_v^2}$

\par Hence our proof is complete.
\qed

\bt \label{ratk}
If $G$ is a connected graph of order $n$ then $\frac{\r_1}{\r_2}\le \frac{n-1}{n-2}, $ equality holds if and only if $G=K_n.$
\et
\pf
If $\x$ be the normalised perron vector of $\D(G)$ with $x_v={\ds \min_{i\in V(G)}x_i},$ then
\be
\label{661} x_v\le \frac{1}{\sqrt n}
\ee

Therefore by \cref{ratio} we have
\be
\label{662} \frac{\r_1}{\r_2}& \le \frac{1-x_v^2}{1-2x_v^2} \\
\label{663} & \le \frac{n-1}{n-2}
\ee
Now equality in \cref{661} holds if and only if $\x=\frac{\1}{\sqrt n}.$ Also by \cref{ratio} equality in \cref{662} holds if and only if
$
x_u=\frac{d_{uv}}{\sqrt{(\r_1-\r_2)(2\r_1-\r_2)}} \quad \forall \, u \neq v
$

Thus equality in \cref{663} holds if and only if

\beq
x_u &=\frac{d_{uv}}{\sqrt{(\r_1-\r_2)(2\r_1-\r_2)}} =\frac{1}{\sqrt n} \qquad \forall \, u \neq v \\
\text{ i.e. } d_{uv}^2 &=\frac{(\r_1-\r_2)(2\r_1-\r_2)}{n} \qquad \forall \, u \neq v \\
\eeq
But $G$ is connected, therefore we must have
$d_{uv}=1 \, \forall \, u\ne v$ and thus $Tr(v)=n-1.$
Again $x_u=\frac{1}{\sqrt n} \, \forall \, u \ne v $ implies $x_v=\frac{1}{\sqrt n}.$ Therefore $\x=\frac{\1}{\sqrt n}.$
Thus $G$ is transmission regular i.e. all the row sums of $D(G)$ are equal. Hence $Tr(v)=n-1$ implies $G=K_n.$
\qed
\vs

\cref{ratio} can also be expressed in slightly different form as follows.

\bl \label{diff}
If $G$ be a connected graph and $\x$ be the normalized perron vector of $\D(G)$ then for any $v \in V(G)$,
\beq
\r_1-\r_2\le \frac{\r_2 x_v^2}{1-2x_v^2}
\eeq
equality holds if and only if $x_u=\frac{d_{uv}}{\sqrt{(\r_1-\r_2)(2\r_1-\r_2)}} \quad \forall \, u\neq v$
\el

Using \cref{diff} and proceeding as in \cref{ratk} the following result can easily be established.

\bt
If $G$ is a connected graph of order $n$ then
$\r_1-\r_2\le \frac{\r_2}{n-2},$ with equality if and only if $G=K_n.$
\et

%%%%%%%%%%%%%%%%%%%%%%%%%%%%%
\section{Smallest five distance Pareto eigenvalues}\label{upto 5th}
In this section, we provide all possible values of the smallest five distance Pareto eigenvalues of a connected graph.

\bt \label{1-3} For any connected graph $G$ with at least $3$ vertices, $0, 1$ and $2$ are the smallest three distance Pareto eigenvalues of $G.$
\et
\pf If $G$ is the complete graph, then the result follows from \cref{k}. Suppose $G$ is not complete. Then 0 being the only $1\times 1$ principal sub-matrix of $\D(G)$ is the smallest Pareto eigenvalue of $\D(G).$ Again as $G$ is connected, $A=\bpm
0 & 1 \\
1 &0
\epm $
is a $2\times 2$ principal sub-matrix of $\D(G)$ and any other $2\times 2$ sub-matrix of $\D(G)$ dominates $A$. Hence $\mu_2(G)=1.$ Now as $G$ has at least 3 vertices and $G$ is not complete, therefore $B=\bpm
0 & 2 \\
2 &0
\epm $
is a $2\times 2$ principal sub-matrix of $\D(G)$ and any other $2\times 2$ principal sub-matrix of $\D(G)$ other than $A$ and $B$ dominates both of them. Further $\r(B)=2$ and any principal sub-matrix of $\D(G)$ of order 3 or higher has minimum row sum 2 and hence spectral radius at least 2. Therefore $\mu_3(G)=2.$

\bt \label{4}
The fourth smallest distance Pareto eigenvalue of a connected non complete graph is $1+\sqrt{3}.$
\et
\pf If $G$ is a connected non complete graph of order $n$, then $n\geq 3$ and therefore from \cref{1-3} we see that 0, 1, 2 are the smallest three distance Pareto eigenvalues of $G.$
Now let $A_1=\J_3-I_3,$ $A_2=\D(P_3).$ Then\\
\beq \r(A_1)=2 \hs \text{and} \hs \r(A_2)=1+\sqrt{3.}\eeq

Now it can be observed from diam$(G)\geq 2,$ therefore any principal sub-matrix of $\D(G)$ of order 3 or higher other than $A_1,A_2$ always dominates either $A_1$ or $A_2.$ Besides any principal sub-matrix of $\D(G)$ of order 4 or higher dominating $A_1$ has minimum row sum 3 and therefore has spectral radius at least 3. Also if diam$(G)\geq 3,$ then $3\in\Pi(G).$ On the other hand $\r(A_2)=1+\sqrt{3}<3.$
\par Again as diam$(G)\geq 2,$ $\D(P_3)=A_2$ is always a principal sub-matrix of $\D(G).$ \\
Hence $\mu_4(G)=1+\sqrt{3}.$
\qed

{\bf Note:} From \cref{k,4}, we observe that a connected graph $G$ is complete if and only if $1+\sqrt{3}\notin \Pi (G).$ Also among all connected graphs of given order $n,$ $K_n$ is the only graph with all integral distance Pareto eigenvalues.

\bt \label{51} If $G$ is a non complete graph with at least $4$ vertices, then
$ \mu_5(G)\geq 3. $ \\
The equality holds if and only if $\omega(G)\geq 4$ or $diam(G)\geq 3.$
\et
\pf Let $A_0=3(\J_2-I_2),$ $A_1=\J_3-I_3,$ $A_2=\D(P_3),$ $A_3=\left(\begin{array}{rrr}
0 & 2 & 2 \\
2 & 0 & 1 \\
2 & 1 & 0
\end{array}\right),$ \\
$A_4=2(\J_3-I_3),$ $B_1=\J_4-I_4,$ $B_2=\D(K_4-e),$ and $B_3=\D(C_4).$

From \cref{4} we have $\mu_4(G)=\r(A_2).$ \\

Now as diam$(G)\geq 2,$ any principal sub-matrix of $\D(G)$ of order 3 or higher other than $A_0,A_1,A_2,A_3,A_4,B_1,B_2,B_3$ dominates at least one of $A_i$ or $B_j$ for $i=1,2,3,4$ and $j=1,2,3.$ Besides any principal sub-matrix of $\D(G)$ of order 4 or higher dominating $A_1$ is either $B_1$ or it dominates $B_1.$ Similarly any principal sub-matrix of $\D(G)$ of order 4 or higher dominating $A_2$ is either $B_2$ or it dominates $B_2$ and hence
dominates $B_1.$
\par Again $\min \{\r(A_0), \r(A_3), \r(A_4), \r(B_1), \r(B_2), \r(B_3) \}=3,$ and equality occurs for $A_0$ and  $B_1.$
\par Hence $\mu_5(G)\geq 3.$ The equality holds if and only if $A_0$ or $B_1$ is a principal sub-matrix of $\D(G)$ which is the case if and only if either $K_4$ or $P_4$ is a induced connected subgraph of $G,$ i.e. if and only if $\omega(G)\geq 4$ or diam$(G)\geq 3.$
\qed

\bt \label{52} If $G$ is a connected graph of order $n\geq 4,$  diameter $2$ and  $\omega(G)\leq 3,$ then
\beq
\mu_5(G)=
\begin{cases}
\frac{1+\sqrt{33}}{2}& \text{ if } C_5 \hs or \hs S_4^+ \text{ is an induced subgraph of G,}\\
\frac{3+\sqrt{17}}{2} &\text{ if neither } C_5 \text{ nor }  S_4^+ \text{ is an induced subgraph of G but } K_4-e \text{ is,}\\
4 &\text{ otherwise.}
\end{cases}
\eeq

\et
\pf If we take $A_0, A_1, A_2, A_3, B_1, B_2, B_3$ as in the \cref{51}, then from the proof of \cref{51} it is clear that
\beq
\mu_5(G) &\geq \min \{\r(A_3), \r(A_4),\r(B_2), \r(B_3) \}\\
&= \r(A_3) \\
&= \frac{1+\sqrt{33}}{2},
\eeq

where equality holds if and only if $A_3$ is a principal sub-matrix of $\D(G),$ i.e. if and only if $C_5$ or $S_4^+$ is an induced subgraph of G.
Now if $\mu_5(G) < \frac{1+\sqrt{33}}{2},$ then neither $C_5$ nor $S_4^+$ is an induced subgraph of $G$ and therefore

\beq
\mu_5(G) &\geq \min \{\r(A_4),\r(B_2), \r(B_3) \}\\
&= \r(B_2) \\
&= \frac{3+\sqrt{17}}{2},
\eeq

where equality holds if and only if $B_2$ is a principal sub-matrix of $\D(G),$ i.e. if and only if $K_4-e$ is an induced subgraph of G.

Finally if $\mu_5(G) < \frac{3+\sqrt{17}}{2},$ then none of $K_5, P_4, S_4^+, K_4-e$ is an induced subgraph of $G$ and therefore either $C_4$ or $S_4$ must be an induced subgraph of $G.$ Now if $C_4$ is an induced subgraph of $G,$ then $B_3$ is a principal sub-matrix of $\D(G)$ and in the other case $A_4$ is a principal sub-matrix of $\D(G).$ In the either case, we have
\beq \r(B_3)=\r(A_4)=4.
\eeq
Hence $\mu_5(G)=4.$
\qed

\bc If $G$ is a non complete connected graph with at least $4$ vertices, then
\beq 3\leq\mu_5(G)\leq 4. \eeq

Furthermore, the left hand equality holds if and only if \mbox{diam}$(G)\geq 3$ or $\omega(G)\leq 4$ and the right hand equality holds if and only if $G$ does not have $K_4,$ $P_4,$ $C_5,$  $K_4-e,$ and $S_4^+$ as induced subgraph.

\ec

\bc \label{5t} If $T$ is a tree with at least $4$ vertices, then $\mu_5(T)$= $4$ or $3$ according as $T$ is a star or not.\ec
\pf For any tree $T,$ we have $\omega(T)=2.$ Therefore if diam$(T)\geq3,$ then by \cref{51}, $\mu_5(T)= 3.$\\

Now if diam$(T)=2,$ T must be a star and then by \cref{52}, $\mu_5(T)= 4.$ \qed

%\bc Among all trees with at least $4$ vertices, 5th smallest distance Pareto eigenvalue is maximum for the star.\ec

\bc If $n\geq 4,$ then
\beq \mu_5(C_n)=
\begin{cases}
4 \text{ if } n=4,\\
\frac{1+\sqrt{33}}{2} \text{ if } n=5,\\
 3 \text{ otherwise. }
\end{cases}
\eeq
\ec

\bc For any positive integers $m,n$ with $m+n\geq 4,$ $\mu_5(K_{m,n})=4.$
\ec
%%%%%%%%%%%%%%%%%%%%%%%%%%%%%%%%%%%%%%%%%%%%%%%%%%%%%%%%
\section{6th smallest distance Pareto eigenvalue}\label{6th}

From \cref{51} and \cref{52}, we see that for a connected graph $G$ possible values of $\mu_5(G)$ are $3,4, \frac{3+\sqrt{17}}{2}$ and $\frac{1+\sqrt{33}}{2}.$ In this section we consider all those four cases and find all possible values of the sixth smallest distance Pareto eigenvalue of a connected graph.

\bt \label{61} If $G$ is a connected graph with at least $5$ vertices and $\mu_5(G)=4,$ then
\beq \mu_6(G)=
\begin{cases}
5 \mbox{  if  } G=K_6,\\
 2+\sqrt{7}  \mbox{ otherwise. }
\end{cases}
\eeq
\et

\pf First suppose that $G$ is a complete graph, then by \cref{k}, $G=K_6$ is the only complete graph with $\mu_5(G)=4$ and in this case $\mu_6(G)=5.$
\par Now if $G$ is not a complete graph then by \cref{51} and \cref{52}, $\mu_5(G)=4$ implies  diam$(G)=2$ and $\omega(G)\leq 3.$ \\

So any principal sub-matrix of $D(G)$ of order at most 3 will have row (column) sum at most 4.
\par Therefore
\beq
 \mu_6(G)=\min_A \r(A),
\eeq
where the minimum is over all principal sub-matrix of $D(G)$ of order 4 or higher with $\r(A)>4.$\\

Now as $\omega(G)\leq 3,$ $\J_k-I_k$ cannot be a principal sub-matrix of $\D(G)$ for $k \geq 5.$

{\bf Claim:} $S_4$ must be an induced subgraph of $G.$\\
$\mu_5(G)=4 $ implies that $G$ does not have $K_4, K_4-e, P_4$ and $S_4^+$ as an induced subgraph. Therefore only possible induced connected subgraphs of order 4 are $C_4$ and $S_4.$ If $S_4$ is an induced subgraph of $G$ then we are done. \\
Otherwise let $H=C_4$ be an induced subgraph of $G.$  Since $G$ has at least 5 vertices, we can choose vertex $w \in V(G)-V(H).$ Again as diam$(G)=2,$ $w$ must be adjacent to at least two vertices of $H.$ But as $S_4^+$ is not an induced subgraph of $G,$  $w$ cannot be adjacent to two adjacent vertices in $H.$ Again for the same reason $w$ cannot be adjacent to more than 2 vertices in $H.$ Hence $w$ must be adjacent to exactly two vertices in $H$ which are not adjacent in $H.$ Thus $S_4$ must be an induced subgraph of $G$ and thereby the claim is established. \\

Now as diam$(G)=2$ and $\omega(G)\leq 3,$ it is obvious that any principal sub-matrix of $\D(G)$ of order 4 or higher other than $\D(S_4)$ always dominates $\D(S_4).$

Hence $\mu_6(G)=\r(D(S_4))=2+\sqrt{7}.$
\qed

\bt \label{62} If $G$ is a connected graph with at least $5$ vertices and $\mu_5(G)=\frac{3+\sqrt{17}}{2},$ then  $\mu_6(G)=4.$ \et
\pf From \cref{52},  $\mu_5(G)=\frac{3+\sqrt{17}}{2}$ implies that $K_4, P_4, C_4$ and $S_4^+$ are not an induced subgraph of $G$ but $K_4-e$ is. As $K_4$ is not induced subgraph of $G$ any principal sub-matrix of $\D(G)$ of order 5 or higher dominates $\J_5-I_5$ and has spectral radius greater than 4.
Besides diam$(G)=2$ implies that spectral radius of any principal sub-matrix of $\D(G)$ of order 1 or 2 is at most 2.
\par Therefore
\beq
\mu_6(G)=\min_A \r(A),
\eeq

where the minimum is over all principal sub-matrix $A$ of $\D(G)$ of order 3 or higher with $\r(A)>\frac{3+\sqrt{17}}{2}.$
\par Now let $H=K_4-e$ is an induced subgraph of $G$ then as $G$ has at least 5 vertices, we can choose vertex $w$ of $G$ such that $w\notin V(H).$ Again diam$(G)=2$ implies that $w$ must be adjacent to at least one vertex of $V(H).$

\par If $w$ is adjacent to two vertices of $H$ both of degree 3, then $S_4$ is an induced subgraph of $G$ and therefore $M=2(\J_3-I_3)$ is a principal sub-matrix of $\D(G).$ Again if $w$ is adjacent to three vertices of $H$ of which two are of degree 2 and the third is of degree 3 then $C_4$ is an induced subgraph of $G$ and therefore $N=\D(C_4)$ is a principal sub-matrix of $\D(G).$\\ Now as $G$ cannot have any of  $S_4^+,K_4, P_4$ as an induced subgraph, therefore it can be easily observed that there cannot be any other possibilities for $w.$
\par But $\r(M)=\r(N)=4.$  Hence $\mu_6(G)=4.$
\qed

\bt \label{63} If $G$ is a connected graph with at least $5$ vertices and $\mu_5(G)=3,$ then
\beq
\mu_6(G)=
\begin{cases}
\frac{1+\sqrt{33}}{2}  &\mbox{  if  } C_5 \text{ or } S_4^+ \text{ is an induced subgraph of G,}\\
\frac{3+\sqrt{17}}{2}  &\mbox{  if  } C_5  \text{ and } S_4^+ \text{ are not induced subgraph of G but } K_4-e \text{ is,}\\
4 \hs &\mbox{  if  } C_5,S_4^+,K_4-e \text{ are not induced subgraph of G but at least } \\
&\text{ one of }  K_5, C_6, C_4, S_4, P_5 \text{ is,}\\
\r(\D(P_4)) &\text{ otherwise.}
\end{cases}
\eeq
\et
\pf From \cref{51}, $\mu_5(G)=3$ implies $\text{diam}(G)=3$ or $\omega(G)\geq 4.$ Proceeding as in the \cref{62} we can show that $\mu_6(G)=\frac{1+\sqrt{33}}{2}$ if $C_5$ or $S_4^+$ is an induced subgraph of G and $\mu_6(G)=\frac{3+\sqrt{17}}{2}$ if $C_5$ and $S_4^+$ are not induced subgraph of $G$ but $K_4-e$ is.
\par Now suppose $C_5, K_4-e$ and $S_4^+$ are not induced subgraph of $G.$ Then as diam$(G)\geq 3,$ or $\omega(G)\geq 4$ or both, therefore $K_4$ or $C_4$ or $S_4$ or $P_4$ must be induced subgraph of $G.$ Thus at least one of $\D(K_4), B_3,A_4, \D(P_4)$ must be a principal sub-matrix of $\D(G),$ where $A_4$ and $B_3$ are as defined in \cref{51}. But $\r(\D(K_4))=3$ and so we can ignore it.

Also $\r(A_4)=\r(B_3)=4<\r(\D(P_4)).$ Besides if $\omega(G)\geq 5$ then $\D(K_5)$ is a principal sub-matrix of $\D(G)$ with $\r(\D(K_5))=4.$ Again $A_4$ is a principal sub-matrix of $\D(G)$ if $C_6$ is a induced subgraph of $G.$ Also as $\r(\D(K_5-e))>\r(\D(P_4))>4$ and $\r(\D(K_6))=5\r(\D(P_4)),$ so for any sub-matrix of $\D(G)$ of order 5 or higher other than $\D(K_5)$ spectral radius is more than $\r(\D(P_4)).$ Hence if at least one of $K_5, C_6, C_4, S_4, P_5$ is an induced subgraph of $G,$ then $\mu_6(G)=4$ and otherwise $\mu_6(G)=\r(\D(P_4)).$
\qed

\bc If $T$ is a tree with $n\geq 5$ vertices, then
\beq \mu_6(T)=
\begin{cases}
2+\sqrt{7}  \mbox{  if  } T=S_n,\\
4  \mbox{  otherwise.}
\end{cases}
\eeq
\ec
\pf If diam$(T)=2,$ then $T$ must be a star, therefore from \cref{5t}, $\mu_5(T)=4$ and thus by \cref{61} $\mu_6(T)=2+\sqrt{7}.$\\
\par Now if diam$(T)= 3,$ then $T$ must have $S_4$ as induced subgraph as $n\geq 5.$ Again if diam$(T)\geq 4,$ then $T$ must have $P_4$ as induced subgraph. Thus in either case by  \cref{63} we have $\mu_6(T)=4.$
\qed

\bc Among all trees with at least 5 vertices, 6th smallest distance Pareto eigenvalue is maximum for the star graph.
\ec

\bt \label{64} If $G$ is a connected graph with at least 5 vertices and $\mu_5(G)=\frac{1+\sqrt{33}}{2},$ then
\beq \mu_6(G)=
\begin{cases}
\frac{3+\sqrt{37}}{2}  \mbox{  if  } G=C_5,\\
\gamma  \mbox{  if  } G=C_3*C_3,\\
\frac{3+\sqrt{17}}{2} \mbox{  if  } K_4-e \text{ is an induced subgraph of } G,\\
4 \text{ otherwise.}
\end{cases}
\eeq
where $\gamma$ is the largest root of $x^3-x^2-11x-7=0$
\et

\pf As before, it can be easily shown that if $K_4-e$ is an induced subgraph of $G$ then $\mu_6(G)=\frac{3+\sqrt{17}}{2}.$\\
\par Now suppose that $K_4-e$ is not  an induced subgraph of $G.$ As $\mu_5(G)=\frac{1+\sqrt{33}}{2},$ from \cref{51} and \cref{52} $G$ does not have $P_4$ or $K_4$ as induced subgraph but has $C_5$ or $S_4^+$ as induced subgraph.
\par If $H=C_5$ is an induced subgraph of $G$ then for $n=5$, $G=C_5$ and $\mu_6(C_5)$ is the spectral radius of any $4\times 4$ sub-matrix of $\D(C_5)$ i.e. $\mu_6(G)=\frac{3+\sqrt{37}}{2}.$ Again if
$n\geq 6,$ then as diam$(G)=2$ and $K_4$ is not an induced subgraph of $G,$ therefore any vertex $w\in V(G)-V(H)$ of $G$ must be adjacent to at least two non adjacent vertices of $H.$ Thus at least one of $C_4$ and $S_4$ must be an induced subgraph of $G.$ Hence $\mu_6(G)=4.$

\par Again if $H=S_4^+$ is an induced subgraph of $G,$ we take $w\in V(G)- V(H).$ Now as diam$(G)=2,$ $w$ must be adjacent to at least one vertex of $H.$ Besides as $K_4-e$ is not induced subgraph of $G,$ $w$ cannot be adjacent to $3$ or more vertices of $H,$ also for the same reason $w$ cannot be adjacent to two vertices in the triangle in $H.$

 If $w$ is adjacent to a single vertex $u\in H,$ then diam$(G)=2$ implies that $d_H(u)=3.$ Thus $S_4$ is an induced subgraph of $G$. Again if $w$ is adjacent to exactly two vertices $u,v\in H,$ with $\{d_H(u), d_H(v) \}=\{1, 2 \}$ then $C_4$ is an induced subgraph of $G.$ In either case we get $\mu_6(G)=4.$
 \par Now if $G$ does not have $C_4$ or $S_4$ as induced subgraph then we are left with only one possibility i.e. every $w\in V(G)- V(H)$ is adjacent  to exactly two vertices $u,v\in H,$ with $\{d_H(u), d_H(v) \}=\{1, 3 \}.$ But in this situation we must have  $|V(G)- V(H)|=1$ as $K_4-e$ is not an induced subgraph of $G.$ Which implies that $G=C_3*C_3.$ It can be directly verified that $\mu_6(C_3*C_3)$ is the spectral radius of the matrix $\left(\begin{array}{rrrr}
0 & 1 & 2 & 2 \\
1 & 0 & 1 & 1 \\
2 & 1 & 0 & 1 \\
2 & 1 & 1 & 0
\end{array}\right). $\\
Therefore $\mu_6(G)$ is the largest root of $x^3-x^2-11x-7=0.$

\par Combining all the above situations, we get our required result.
\qed

\bc If $n\geq 5,$ then
\beq \mu_6(C_n)=
\begin{cases}
\frac{3+\sqrt{37}}{2} \text{ if } n=5,\\
\gamma \text{ if } n=7,\\
4 \text{ otherwise. }
\end{cases}
\eeq
where $\gamma$ is the largest root of $x^3-x^2-11x-7=0$
\ec

\bc For any positive integers $m,n$ with $m+n\geq 5,$ $\mu_6(K_{m,n})=2+\sqrt{7}.$
\ec

%%%%%%%%%%%

%%%%%%%%%%%%%%%%%%%%%%%%%%%%%%%%%%%%%%%
%%%%%%%%%%%%%%%%%%%%%%%%%%%%%%%%%%%%%%
\comment{
From \cref{r2,ratio,thmsn,r1r2} we get.

Now \Vref{r1r2,ratio} and \crefrange{r1r2}{ratio}.

Also \cpageref{r1r2,ratio} }
%%%%%%%%%%%%%%%%%%%%%%%%%%%%%%%%

\end{document}